\documentclass[12pt]{amsart}
\usepackage{amscd}
\include{epsf}

\usepackage{fancyhdr}
\usepackage[all]{xy}
\usepackage{amssymb}
\newtheorem{Theorem}{Theorem}

\newtheorem{Definition}[Theorem]{Definition}

\newcommand{\QED}{{\hfill$\Box$\medskip}}

\def \0{\lambda_{0}}

\def\Fs{{\mathcal F}}

\def\is{{\rm Iso}}

\def\co{\colon\thinspace}

\usepackage{graphicx, amssymb, amsfonts}

\usepackage{epsfig}

\begin{document}

\title{The volume flux group and nonpositive curvature}

\date{\today}
\thanks{The author is supported by the Deutsche Forschungsgemeinschaft (DFG) project `Asymptotic Invariants of Manifolds'.}
\address{Mathematisches Institut LMU, Theresienstrasse 39, M\"unchen 80333 Germany.}
\address{CIMAT, CP: 36240, Guanajuato, Gto, M\'exico.  }
\email{p.suarez-serrato@cantab.net}


\maketitle

\begin{center}
by  Pablo Su\'arez-Serrato \\ [.5cm]
 \end{center}
\begin{abstract}
We show that every closed nonpositively curved manifold with non-trivial volume flux group has zero minimal volume, and admits a finite covering with circle actions whose orbits are homologically essential. This proves a conjecture of Kedra--Kotschick--Morita for this class of manifolds.\\ [.5cm]
\end{abstract}


 Let $M$ be a closed smooth manifold and $\mu$ a volume form on $M$. Denote by ${\rm Diff}^\mu$ the group of $\mu$--preserving diffeomorphisms of $M$, and by ${\rm Diff}^\mu_0$ its identity component. The $\mu$--flux homomorphism ${\rm Flux}_\mu$, from the universal covering $\widetilde{\rm Diff}{}^\mu_0$ to the ${(n-1)}$--cohomology group $H^{n-1}(M;{\bf R})$, is defined by the formula
 \[ {\rm Flux}_\mu ([\varphi_t]) = \int_0^1 [i_{\dot{\varphi}_t} \mu] dt.\]
  It induces a homomorphism
  \[ {\rm Flux}_\mu \colon \pi_ 1 ({\rm Diff}^\mu_0) \rightarrow H^{n-1} (M; {\bf R})\]
   whose image is the volume flux group $\Gamma_\mu \subset H^{n-1}(M; {\bf R})$. The $\mu$--flux homomorphism descends to a homomorphism
   \[ {\rm Flux}_\mu \colon {\rm Diff}^\mu_0 \rightarrow H^{n-1} (M; {\bf R})/\Gamma_\mu.\]

For a closed connected smooth Riemannian manifold $(M,g)$, let  ${\rm Vol}(M,g)$ denote the volume of $g$ and let $K_{g}$ be its sectional curvature. We define the minimal volume of $M$ following J Cheeger and M Gromov \cite{CG}:
\[ {\rm MinVol }(M):=\inf\limits_{g} \{ {\rm Vol}(M,g)\; : \; |K_{g}| \leq 1 \} \]
The minimal volume is a very sensitive invariant, it was first observed by L Bessi\`eres \cite{LBes} that its value may depend on the differentiable structure of $(M,g)$. Indeed, D Kotschick \cite{K} has shown that even the vanishing of ${\rm MinVol }(M)$ can detect changes in the smooth structure of $M$.

The investigation of relationships between the volume flux group of $M$ and various invariants which bound ${\rm MinVol }(M)$ from below was set in motion by J Kedra, D Kotschick and S Morita \cite{KKM}. They put forward the idea that if a closed manifold has non-trivial volume flux group then its minimal volume should vanish. The aim of this note is to verify that statemtent for closed manifolds which carry a metric of nonpositive sectional curvature.

The volume flux group $\Gamma_{\mu}$ is independent of the form $\mu$---see section $3$ of \cite{KKM}---so it can be considered as an invariant of the manifold $M$ itself.

\begin{Theorem}
Every closed nonpositively curved manifold with non-trivial volume flux group has zero minimal volume.
\end{Theorem}

The technique we will use to show that the minimal volume vanishes is an $\Fs$--structure, which was introduced by Gromov as a generalisation of an $S^1$--action.

\begin{Definition}{ An $\Fs$--structure on a closed manifold $M$ is given by the following conditions.
\begin{enumerate}
\item{ A finite open cover $\{ U_1, ..., U_{N} \} $}
\item{ $\pi_{i}\co \widetilde{U}_{i}\rightarrow U_{i}$ a finite Galois covering with group of deck transformations $\Gamma_{i}$, $1\leq i \leq N$}

\item{ A smooth torus action with finite kernel of the $k_{i}$-dimensional torus, \\ $\phi_{i}\co T^{k_{i}}\rightarrow {\rm{Diff}}(\widetilde{U}_{i})$, $1\leq i \leq N$}

\item{ A homomorphism $\Psi_{i}\co\Gamma_{i}\rightarrow {\rm{Aut}}(T^{k_{i}})$ such that
\[ \gamma(\phi_{i}(t)(x))=\phi_{i}(\Psi_{i}(\gamma)(t))(\gamma x) \]
for all $\gamma \in \Gamma_{i}$, $t \in T^{k_{i}}$ and $x \in \widetilde{U}_{i}$ }

\item{ For any finite sub-collection $\{ U_{i_{1}}, ..., U_{i_{l}} \} $ such that  $U_{i_{1}\ldots i_{l}}:=U_{i_{1}}\cap \ldots \cap U_{i_{l}}\neq\emptyset$ the following compatibility condition holds: let $\widetilde{U}_{i_{1}\ldots i_{l}}$ be the set of points $(x_{i_{1}}, \ldots , x_{i_{l}})\in \widetilde{U}_{i_{1}}\times \ldots \times \widetilde{U}_{i_{l}}$ such that $\pi_{i_{1}}(x_{i_{1}})=\ldots = \pi_{i_{l}}(x_{i_{l}})$. The set $\widetilde{U}_{i_{1}\ldots i_{l}}$ covers $\pi_{i_{j}}^{-1}(U_{i_{1}\ldots i_{l}}) \subset \widetilde{U}_{i_{j}}$ for all $1\leq j \leq l$, then we require that $\phi_{i_{j}}$ leaves $\pi_{i_{j}}^{-1}(U_{i_{1}\ldots i_{l}})$ invariant and it lifts to an action on $\widetilde{U}_{i_{1}\ldots i_{l}}$ such that all lifted actions commute }
\end{enumerate}}\end{Definition}

An $\Fs$--structure is said to be {\em pure} if all the orbits of all actions at a point, for every point have the same dimension. We will say an $\Fs$--structure is {\em polarised} if the smooth torus action $\phi_{i}$ above are fixed point free for every $U_{i}$. The existence of a polarised $\Fs$--structure on a manifold $M$ implies the minimal volume ${\rm MinVol}(M)$ is zero by the main result of Cheeger and Gromov \cite{CG}, the interested reader is invited to consult the illuminating examples found there as well.

The attentive reader will notice that the definition of an $\Fs$--structure above is different from the sophisticated one found in \cite{CG}. Despite this, it is sufficiently practical to be implemented and also satisfies the properties needed in the proof that a polarised $\Fs$--structure forces the minimal volume to vanish, which can be consulted in \cite{CG}.

\medskip

{\it Proof of Theorem 1.}  Let $M$ be a compact nonpositively curved manifold whose volume flux group is not trivial. The fundamental group $\pi_1(M)$ has non-trivial centre $Z$, a proof can be found in \cite[Theorem 15]{KKM} and compared with work of A Fathi \cite[Proposition 5.1]{Fat} and the various references and attributions contained therein. In this case $M$ admits a finite covering space $M^{\ast}$ diffeomorphic to $T^{k} \times N$, where $T^{k}$ is a flat torus of dimension $k$ and $N$ is a compact nonpositively curved manifold as was shown by P Eberlein \cite{Eb1}.

Even though it may seem that this already implies ${\rm MinVol}(M)=0$, for completeness will now show how to construct a pure polarised $\Fs$--structure on $M$; since the torus $T^{k}$ splits off $M^{\ast}$ smoothly and the action of $T^{k}$ on itself as a factor of $M^{\ast}$ is compatible with the covering transfomation in the required sense. This will also provide an example of a detailed construction of a pure polarised  $\Fs$--structure.

Represent $M$ as $H / \Gamma $, where $H$ is simply connected and $\Gamma $ is a properly discontinuous group of isometries of $H$ which acts freely. The space $H$ decomposes into $H_1 \times H_2$, where $H_1= {\bf R}^{k}$ is a Euclidean space of dimension $k= {\rm rank}(Z)$. 

Notice that every element $\gamma$ of $\Gamma $ is of the form $\gamma = \gamma_1 \times \gamma_2 \in \is(H_1)\times \is(H_2)$, here $\is(H_{i})$ denotes the group of isometries of $H_{i}$ (see Lemma 1 in \cite{Eb1}, and the subsequent discussion). Let $p_{i}\co \Gamma \to \is(H_{i})$ denote the projection homomorphisms, then $\Gamma_1 = p_1(\Gamma )$ acts by translations on $H_1$ and $\Gamma_2=p_2(\Gamma)$ is a discrete subgroup of $\is (H_2)$. So $Z \subset \Gamma_1$ and $H_1 / Z$ is a compact flat torus $T^{k}$. The projection $p \co H_1 \to T^{k}$ allows us to define $\rho \co \Gamma_2 \to T^{k}$ by setting $\rho(p_2 \gamma )= p(p_1 \gamma )$. The function $\rho $ is well defined since $\ker (p)= \ker (p_2)= Z$. The centre $Z$ can also be thought of as a set of vectors in $H_1$, as $Z \subset \Gamma_1$. In this way $Z$ acts on $H_1$ by translations.

Recall that $\Gamma_2$ has trivial centre and that $M$ is isometric to $(T^{k}\times H_2) / \Gamma_2$ \cite{Eb1}, here $\Gamma_2$ acts on $T^{k}\times H_2$ by $\psi (s,h)=(\rho(\psi)s, \psi(h))$ with $(s,h)$ in $T^{k}\times H_2$ and $\psi$ in $\Gamma_2$. This can be read from the diagram

\[ \xymatrix{H=H_1 \times H_2 \ar[r]^{p \times Id} \ar[d]  & T^{k}\times H_2 \ar[d]^{q} \\
            M=H / \Gamma \ar[r]^{F} &  (T^{k}\times H_2) / \Gamma_2 } \]
where $F$ is defined so that the diagram commutes.

In $\Gamma_2$ there exists a finite index subgroup $\Gamma_0$---denoted by $\Gamma_2^{
\ast \ast}$ in \cite{Eb1}---which makes the following diagram commute.

\[ \xymatrix{ H \ar[r] \ar[dr] & T^{k}  \times ( H_2 / \Gamma_0 ):= M^{\ast} \ar[d] \\
                     & (T^{k} \times H_2)/ \Gamma_2 = M } \]
Define $q^{\ast}\co T^{k}  \times ( H_2 / \Gamma_0 ) = M^{\ast} \to M$ as in the previous diagram and denote by $\Gamma^{\ast}$ the group of deck transformations of $q^{\ast}$ seen as a covering map. Notice that $\Gamma^{\ast} \subset \Gamma_2$.

The function $\rho$ is defined on all of $\Gamma_2$, so we can also consider the restriction of $\rho$ to $\Gamma^{\ast}$.

 We are now in a position to verify that this construction gives $M$ a pure polarised $\Fs$--structure, and hence ${\rm MinVol}(M)=0$ as claimed. Let us check that every condition which guarantees the existence of a polarised $\Fs$--structure is met.
\begin{enumerate}
\item Take $U=M$, as an open cover with a single set.
\item  Define $\tilde{U}:=  M^{\ast}=T^{k}\times N$, here $N=H_2 / \Gamma_0$. The quotient map
\[ q^{\ast} \co \tilde{U}=T^{k}\times N \to (T^{k}\times H_2) / \Gamma_2 = M = U \]
is a finite Galois covering with deck transformation group $\Gamma^{\ast}$.
\item The $k$--torus $T^{k}$ acts smoothly and without fixed points on itself as a factor of $T^{k}\times (H_2 / \Gamma_0)$. So we have the action
\[ \phi \co T^{k} \to {\rm Diff}(\tilde{U})= {\rm Diff}(T^{k}\times N) \]
 given by $\phi(t)(s,h)=(s+t,h)$.
 \item We will use the function $\rho $ to define the automorphism $\Psi \co \Gamma^{\ast} \to {\rm Aut}(T^{k})$, set $\Psi(\gamma):=\rho(\gamma)$. Let $t \in T^{k}$ and $x\in \tilde{U}=T^{k}\times N \Rightarrow x=(s,h)\in T^{k}\times N$. Notice that for $s$ and $t$ in $T^{k}$ we have that $\rho(\gamma)(s+t)=(\rho(\gamma)s+\rho(\gamma)t )$ since $\rho(\gamma) \in T^{k}$.

 We plug in this information to obtain the following equalities.

 \begin{eqnarray*}
 \gamma( \phi(t)x) & = & \gamma(\phi(t)(s,h)) \\
 		   & = & \gamma(s+t, h) \\
		   & = & (\rho(\gamma)(s+t), \gamma h) \\
		   & = & (\rho(\gamma)s+\rho(\gamma)t, \gamma h )= \star
 \end{eqnarray*}

\begin{eqnarray*}
\phi(\Psi(\gamma)t)(\gamma x) & = & \phi(\rho(\gamma)t)(\gamma x) \\
			      & = & \phi(\rho(\gamma)t)(\rho(\gamma)s, \gamma h) \\
			      & = & (\rho(\gamma)s+\rho(\gamma)t, \gamma h )= \star
\end{eqnarray*}

Therefore $\gamma(\phi_{i}(t)(x))=\phi_{i}(\Psi_{i}(\gamma)(t))(\gamma x) $ and the condition is satisfied.

\item This condition does not need to be verified, because we only have one covering set.
\end{enumerate}
Since the action of $T^{k}$ on itself as a factor of $ T^{k}  \times ( H_2 / \Gamma_0 )$ is fixed point free, the above construction gives $M$ a pure polarised $\Fs$--structure.\QED 

We will state the contrapositive statement to the theorem because it strengthens Corollary 17 of \cite{KKM}: {\it Let $M$ be a closed nonpositively curved manifold with positive minimal volume. Then the volume flux group $\Gamma_{\mu}$ is trivial for every volume form  $\mu$ on $M$.}

Another noteworthy feature of the covering $M^{\ast} \to M$ is that on $M^{\ast} \cong T^{k}\times N$ elements of the volume flux group which come from rotations on the $T^{k}$ factor are circle actions with homologically essential orbits. It is not yet clear if in the general case a non-trivial volume flux group implies existence of a circle action with homologically essential orbits at least in a multiple cover---compare with Remark 19 in \cite{KKM}---but it is virtually so for every closed nonpositively curved manifold $M$, since it is true for $M^{\ast}$.

{\bf Acknowledgements:} I wish to warmly thank Dieter Kotschick for a number of interesting conversations and for commenting on a previous version of this work.

\end{document}